\documentclass[12pt]{amsart}
\usepackage{amssymb}

\newtheorem{theorem}{Theorem}[section]
\newtheorem{lemma}[theorem]{Lemma}

\newtheorem{corollary}[theorem]{Corollary}
\newtheorem{proposition}[theorem]{Proposition}

\newtheorem{lem-def}[theorem]{Lemma-Definition}

\newcommand{\hooklongrightarrow}{\lhook\joinrel\longrightarrow}
\renewenvironment{proof}{{\bfseries Proof.}}{\qed}
\newcommand{\R}{\mathbb R}

\newcommand{\N}{\mathbb N}
\newcommand{\Z}{\mathbb Z}
\newcommand{\Q}{\mathbb Q}

\newcommand{\F}{\mathbb F}

\def\op{\operatorname}

\def\al{\alpha}

\def\ars#1{\renewcommand\arraystretch{#1}}
\def\as{\op{AS}}

\def\cc{\mathcal{C}}

\def\cfa{\left(\ri\right)_{i\in A}}

\def\comb#1#2{\ars{0.9}\left(\!\!\begin{array}{c}
#1\\#2
\end{array}\!\!\right)\ars{1}}
\def\cl#1{\left[\,#1\,\right]_\mu}

\def\d{\Delta}
\def\defn{\nn{\bf Definition. }}
\def\dep{\op{\mbox{\small depth}}}

\def\dg{\op{deg}}

\def\diso{\lower.4ex\hbox{$\downarrow$}\raise.4ex\hbox{\mbox{\scriptsize
$\wr$}}}

\def\dta{\delta}

\def\e{\medskip}

\def\ep#1{\exp(\Pi i#1)}
\def\ep{\epsilon}

\def\g{\Gamma}
\def\ga{\gamma}

\def\gal{\op{Gal}}

\def\gen#1{\big\langle\, {#1} \,\big\rangle}

\def\gg{\mathcal{G}}

\def\ggm{\mathcal{G}_\mu}

\def\gm{\g_\mu}

\def\gq{\g_\Q}

\def\hra{\hooklongrightarrow}

\def\imp{\ \Longrightarrow\ }

\def\inm{\op{in}_\mu}
\def\inn{\op{in}}

\def\irr{\op{Irr}}
\def\ism{\lower.3ex\hbox{\ars{.08}$\begin{array}{c}\,\to\\\mbox{\tiny $\sim\,$}\end{array}$}}
\def\iso{\ \lower.3ex\hbox{\ars{.08}$\begin{array}{c}\lra\\\mbox{\tiny $\sim\,$}\end{array}$}\ }

\def\kb{\overline{K}}
\def\km{k_\mu}

\def\kp{\op{KP}}

\def\kpi{\op{KP}_\infty}
\def\kpm{\op{KP}(\mu)}

\def\kx{K[x]}

\def\ldp{\op{\mbox{\small lim-depth}}}

\def\lg{l\raise.6ex\hbox to.2em{\hss.\hss}l}

\def\lra{\,\longrightarrow\,}

\def\lx{\operatorname{lex}}

\def\mlv{MacLane--Vaqui\'e\ }

\def\mmu{\mid_\mu}

\def\nn{\noindent}

\def\om{\omega}
\def\oo{\mathcal{O}}
\def\orb{\hbox to  .3em{$\backslash$}\backslash}

\def\ppa{\mathcal{P}_{\al}}

\def\pset{\mathcal{P}}

\def\rhc{\rho_\cc}
\def\ri{\rho_i}

\def\rj{\rho_j}

\def\rni{\rho_{n,i}}

\def\sg{\sigma}

\def\sii{\ \Longleftrightarrow\ }

\def\sni{s_{n,i}}

\def\str{^{\str}}

\def\str{^{\op{str}}}
\def\supp{\op{supp}}

\def\tmn{\ty(\mu,\nu)}

\def\ttt{\mathcal{T}}
\def\ty{\mathbf{t}}

\def\vni{v_{n,i}}

\newcounter{cs}
\stepcounter{cs}
\newcommand{\casos}{\begin{itemize}}
\newcommand{\fcasos}{\end{itemize}\setcounter{cs}{1}}

\newfont{\tit}{cmr12 scaled \magstep3}

\title[Infinite limit-depth]{Valuations with infinite limit-depth}

\makeatletter
\@namedef{subjclassname@2010}{%
  \textup{2010} Mathematics Subject Classification}

\author[Alberich]{Maria Alberich-Carrami$\tilde{\mbox{n}}$ana}
\address{Institut de Rob\`otica i Inform\`atica Industrial (IRI, CSIC-UPC), Institut de Mate\-mà\-tiques de la UPC-BarcelonaTech (IMTech) and Departament de Mate\-m\`a\-tiques, Universitat Polit\`ecnica de Cata\-lunya $\cdot$ BarcelonaTech, Av. Diagonal, 647, E-08028 Barcelona, Catalonia}
\email{Maria.Alberich@upc.edu}
\author[Gu\`ardia]{Jordi Gu\`ardia}
\address{Departament de Matem\`atiques, Escola Polit\`ecnica Superior d'Enginye\-ria de Vilanova i la Geltr\'u, Av. V\'\i ctor Balaguer s/n. E-08800 Vilanova i la Geltr\'u, Catalonia}
\email{jordi.guardia-rubies@upc.edu}
\author[Nart]{Enric Nart}
\author[Ro\'e]{Joaquim Ro\'e}
\address{Departament de Matem\`{a}tiques,         Universitat Aut\`{o}noma de Barcelona,         Edifici C, E-08193 Bellaterra, Barcelona, Catalonia}
\email{nart@mat.uab.cat,\quad jroe@mat.uab.cat}

\thanks{Partially supported by grants  PID2020-116542GB-I00 and PID2019-103849GB-I00 from the Spanish Research Agency, and grant 2017SGR-932 from Generalitat de Catalunya}
\date{\today}

\begin{document}


\begin{abstract}
For a certain field $K$, we construct a valuation-algebraic valuation on the polynomial ring $\kx$, whose underlying Maclane--Vaqui\'e chain consists of an infinite (countable) number of limit augmentations. 
\end{abstract}

\maketitle



\section*{Introduction}
Let $(K,v)$ be a valued field. 
In a pioneering work, Maclane studied the
extensions of the valuation $v$ to the polynomial ring $\kx$ in the case $v$ discrete of rank one \cite{mcla}. 
He proved that all extensions of $v$ to $\kx$ can be obtained as a kind of limit of chains of augmented valuations:
\begin{equation}\label{depthintro}
\mu_0\ \stackrel{\phi_1,\ga_1}\lra\  \mu_1\ \stackrel{\phi_2,\ga_2}\lra\ \cdots
\ \lra\ \mu_{n-1} 
\ \stackrel{\phi_{n},\ga_{n}}\lra\ \mu_{n}\ \lra\ \cdots\ \lra\ \mu
\end{equation}
involving the choice of certain \emph{key polynomials} $\phi_n \in K[x]$ and elements $\ga_n$ belonging to some extension of the value group of $v$.

These chains of valuations contain relevant information on $\mu$ and play a crucial role in the resolution of many  arithmetic-geometric tasks in number fields and function fields of curves  \cite{newapp,gen}.

For valued fields of arbitrary rank, several approaches to this problem were developed by Alexandru-Popescu-Zaharescu \cite{APZ},  Kuhlmann \cite{Kuhl},  
Herrera-Mahboub-Olalla-Spivakovsky \cite{hos,hmos} and Vaqui\'e \cite{Vaq,Vaq3}. 

In this general context, \emph{limit augmentations} and the corresponding \emph{limit key polynomials} appear as a new feature. In the henselian case, limit augmentations are linked with the existence of \emph{defect} in the extension $\mu/v$ \cite{VaqDef}. Thus, they are an obstacle for local uniformization in positive characteristic.

A chain as in (\ref{depthintro}) is said to be a \emph{\mlv chain} if it is  constructed as a mixture of ordinary and limit augmentations, and satisfies certain additional technical condition (see Section \ref{subsecMLV}). In this case, the intermediate valuations $\mu_n$ are essentially unique and contain intrinsic information about the valuation $\mu$ \cite[Thm. 4.7]{MLV}.  
 
In particular, the number of limit augmentations of any \mlv chain of $\mu$ is an intrinsic datum of $\mu$, which is called the \emph{limit-depth} of $\mu$.

In this paper, we exhibit an example of a valuation with an infinite limit-depth,
inspired in a construction by Kuhlmann of infinite towers of Artin-Schreier extensions with defect \cite{KuhlDefect}.

\section{Maclane--Vaqui\'e chains of valuations on $\kx$}\label{secComm}

In this section we recall some well-known results on valuations on a polynomial ring, mainly extracted from the surveys \cite{KP} and \cite{MLV}. 

Let $(K,v)$ be a valued field, with valuation ring $\oo_v$ and residue class field $k$. Let $\g=v(K^*)$ be the value group and denote by $\gq=\g\otimes\Q$ the divisible hull of $\g$. 
 In the sequel, we write $\gq\infty$ instead of $\gq\cup\{\infty\}$.

Consider the set $\ttt$  of all $\gq$-valued extensions of $v$ to the field $K(x)$ of rational functions in one indeterminate. That is, an element $\mu\in\ttt$ is a valuation on $\kx$,
$$
\mu\colon \kx\lra \gq\infty,
$$
such that $\mu_{\mid K}=v$ and $\mu^{-1}(\infty)=\{0\}$.
Let $\gm=\mu(K(x)^*)$ be the value group and $\km$ the residue field.

This set $\ttt$ admits a partial ordering. 
For $\mu,\nu\in \ttt$ we say that $\mu\le\nu$ if $$\mu(f)\le \nu(f), \quad\forall\,f\in\kx.$$
As usual, we write $\mu<\nu$ to indicate that $\mu\le\nu$ and $\mu\ne\nu$.

This poset $\ttt$ has the structure of a tree; that is, all intervals 
$$
(-\infty,\mu\,]:=\left\{\rho\in\ttt\mid \rho\le\mu\right\}
$$
are totally ordered \cite[Thm. 2.4]{MLV}. 

A node $\mu\in\ttt$ is a \emph{leaf} if it  is a maximal element with respect to the ordering $\le$. Otherwise, we say that $\mu$ is an \emph{inner node}.

The leaves of $\ttt$ are the \emph{valuation-algebraic} valuations in Kuhlmann's terminology \cite{Kuhl}.
The inner nodes are the \emph{residually transcendental} valuations, characterized by the fact that the extension $\km/k$ is transcendental. In this case, its transcendence degree is necessarily equal to one \cite{Kuhl}.  

\subsection{Graded algebra and key polynomials}\label{subsecKP}
Take any $\mu\in\ttt$. For all $\alpha\in\g_\mu$, consider the $\oo_v$-modules:
$$
\ppa=\{g\in \kx\mid \mu(g)\ge \alpha\}\supset
\ppa^+=\{g\in \kx\mid \mu(g)> \alpha\}.
$$    
The \emph{graded algebra of $\mu$} is the integral domain:
$$
\ggm=\bigoplus\nolimits_{\alpha\in\g_\mu}\ppa/\ppa^+.
$$

There is an \emph{initial term} mapping $\inm\colon \kx\to \ggm$, given by $\inm0=0$ and  
$$
\inm g= g+\pset_{\mu(g)}^+\quad\mbox{for all nonzero }g\in\kx. 
$$


The following definitions translate properties of the action of  $\mu$ on $\kx$ into algebraic relationships in the graded algebra $\ggm$.\e

\defn Let $g,\,h\in \kx$.

We say that $g$ is \emph{$\mu$-divisible} by $h$, and we write $h\mmu g$, if $\inm h\mid \inm g$ in $\ggm$.

We say that $g$ is $\mu$-irreducible if $\inm g$ is a prime element; that is, the homogeneous principal ideal of $\ggm$ generated by  $\inm g$ is a prime ideal. 

We say that $g$ is $\mu$-minimal if $g\nmid_\mu f$ for all nonzero $f\in \kx$ with $\deg(f)<\deg(g)$.\e

Let us recall a well-known characterization of $\mu$-minimality \cite[Prop. 2.3]{KP}.

\begin{lemma}\label{minimal0}
A polynomial  $g\in \kx\setminus K$ is $\mu$-minimal if and only if $\mu$ acts as follows  on $g$-expansions:
$$
f=\sum\nolimits_{0\le n}a_n g^n,\quad \deg(a_n)<\deg(g)\quad\imp\quad \mu(f)=\min_{0\le n}\left\{\mu\left(a_n g^n\right)\right\}.
$$
\end{lemma}

\defn A  \emph{(Maclane-Vaqui\'e) key polynomial} for $\mu$ is a monic polynomial in $\kx$ which is simultaneously  $\mu$-minimal and $\mu$-irreducible. 
The set of key polynomials for $\mu$ is denoted $\kpm$. \e

All $\phi\in\kpm$ are irreducible in $\kx$.  For all $\phi\in\kpm$ let $\cl{\phi}\subset \kpm$ be the subset of all key polynomials $\varphi\in\kpm$ such that  $\inm \varphi=\inm \phi$.

\begin{lemma}\cite[Thm. 1.15]{Vaq}\label{propertiesTMN}
	Let $\mu<\nu$ be two nodes in $\ttt$. Let $\tmn$ be the set of monic polynomials $\phi\in\kx$ of minimal degree satisfying $\mu(\phi)<\nu(\phi)$.
Then, $\tmn\subset\kpm$ and  $\tmn=\cl{\phi}$ for all $\phi\in\tmn$.

	Moreover,  for all $f\in\kx$, the equality $\mu(f)=\nu(f)$ holds if and only if $\phi\nmid_{\mu}f$.   
\end{lemma}

The existence of key polynomials characterizes the inner nodes of $\ttt$.

\begin{theorem}\label{leaves}
A node $\mu\in\ttt$ is a leaf if and only if $\kpm=\emptyset$.
\end{theorem}

\defn
The \emph{degree} $\deg(\mu)$ of an inner node $\mu\in\ttt$ is defined as the minimal degree of a key polynomial for $\mu$.

\subsection{Depth zero valuations}\label{subsecDepth0}

For all $a\in K$, $\ga\in\gq$, consider the \emph{depth-zero} valuation  $$\mu=\om_{a,\dta}=[v;\,x-a,\ga]\in\ttt,$$
defined in terms of $(x-a)$-expansions as
$$
f=\sum\nolimits_{0\le n}a_n(x-a)^n\imp \mu(f)=\min\{v(a_n)+n\ga\mid 0\le n\}.
$$
Note that $\mu(x-a)=\ga$. 
Clearly, $x-a$ is a key polynomial for $\mu$ of minimal degree and  $\gm=\gen{\g,\ga}$. In particular, $\mu$ is an inner node of $\ttt$ with $\deg(\mu)=1$.

One checks easily that
\begin{equation}\label{balls}
\om_{a,\dta}\le \om_{b,\ep} \ \sii\ v(a-b)\ge\dta \mbox{\, and \,}\ep\ge\dta.
\end{equation}
\subsection{Ordinary augmentation of valuations}\label{subsecOrdAugm}

Let $\mu\in\ttt$ be an inner node.
For all $\phi\in\kpm$ and all $\ga\in\gq$ such that $\mu(\phi)<\ga$, we may construct the \emph{ordinary} augmented valuation  $$\mu'=[\mu;\,\phi,\ga]\in\ttt,$$
defined in terms of $\phi$-expansions as
$$
f=\sum\nolimits_{0\le n}a_n\phi^n,\quad \deg(a_n)<\deg(\phi)\imp \mu'(f)=\min\{\mu(a_n)+n\ga\mid 0\le n\},
$$
Note that $\mu'(\phi)=\ga$, $\mu<\mu'$ and  $\ty(\mu,\mu')=[\phi]_\mu$. 

By \cite[Cor. 7.3]{KP}, $\phi$ is a key polynomial for $\mu'$ of minimal degree. In particular, $\mu'$ is an inner node of $\ttt$ too, with $\deg(\mu')=\deg(\phi)$. 

\subsection{Limit augmentation of valuations}\label{subsecLimAugm}
Consider a totally ordered family of inner nodes of $\ttt$, not containing a maximal element:
$$
\cc=\cfa\subset\ttt.
$$
We assume that $\cc$ is parametrized by a totally ordered set $A$ of indices such that  the mapping
$A\to\cc$ determined by $i\mapsto \ri$
is an isomorphism of totally ordered sets. 

If $\deg(\rho_{i})$ is stable for all sufficiently large $i\in A$, we say that $\cc$ has \emph{stable degree}, and we denote this stable degree by $\dg(\cc)$. 

We say that $f\in\kx$ is \emph{$\cc$-stable} if, for some index $i\in A$, it satisfies
$$\ri(f)=\rj(f), \quad \mbox{ for all }\ j>i.$$

\begin{lemma}\label{stable=unit}
A nonzero $f\in\kx$ is $\cc$-stable if and only if $\,\inn_{\ri}f$ is a unit in $\gg_{\ri}$ for some $i\in A$. 
\end{lemma}

\begin{proof}
Suppose that $\,\inn_{\ri}f$ is a unit in $\gg_{\ri}$ for some $i\in A$. Take any $j>i$ in $A$, and let $\ty(\ri,\rj)=\left[\varphi\right]_{\ri}$. By Lemma \ref{propertiesTMN}, 
$\varphi\in\kp(\ri)$, so that $\inn_{\ri}\varphi$ is a prime element. Hence, $\varphi \nmid_{\ri}f$, and this implies  $\ri(f)=\rj(f)$, again by Lemma \ref{propertiesTMN}.

Conversely, if $f$ is $\cc$-stable, there exists $i_0\in A$ such that $\rho_{i_0}(f)=\ri(f)$ for all $i>i_0$. Hence, $\inn_{\ri}f$ is the image of $\inn_{\rho_{i_0}}f$ under the canonical homomorphism $\gg_{\rho_{i_0}}\to\gg_{\ri}$. By \cite[Cor. 2.6]{MLV}, $\,\inn_{\ri}f$ is a unit in $\gg_{\ri}$.  
\end{proof}\e

We obtain a \emph{stability function} $\rhc$, defined
on the set of all $\cc$-stable polynomials by $$\rhc(f)=\max\{\ri(f)\mid i\in A\}.$$

\defn We say that $\cc$ has a \emph{stable limit} if all polynomials in $\kx$ are $\cc$-stable. In this case, $\rhc$ is a valuation in $\ttt$ and we say that 
$$
\rhc=\lim(\cc)=\lim_{i\in A}\ri.
$$

Suppose that  $\cc$ has no stable limit. Let $\kpi(\cc)$ be the set of all monic $\cc$-unstable polynomials of minimal degree. The elements in $\kpi(\cc)$ are said to be \emph{limit key polynomials} for $\cc$. Since the product of stable polynomials is stable, all limit key polynomials are irreducible in $ \kx$.\e

\defn
We say that $\cc$ is an \emph{essential continuous family} of valuations  if it  has stable degree and it admits limit key polynomials whose degree is greater than $\dg(\cc)$.\e

For all limit key polynomials $\phi\in\kpi\left(\cc\right)$, and all $\ga\in\gq$ such that $\ri(\phi)<\ga$ for all $i\in A$, we may construct the \emph{limit augmented} valuation   $$\mu=[\cc;\,\phi,\ga]\in\ttt$$
defined in terms of $\phi$-expansions as: 
$$ f=\sum\nolimits_{0\le n}a_n\phi^n,\ \deg(a_n)<\deg(\phi)\ \imp \ \mu(f)=\min\{\rhc(a_n)+n\ga\mid 0\le n\}.
$$
Since $\deg(a_n)<\deg(\phi)$, all coefficients $a_n$ are $\cc$-stable. 
Note that $\mu(\phi)=\ga$ and $\ri<\mu$ for all $i\in A$.
By \cite[Cor. 7.13]{KP}, $\phi$ is a key polynomial for $\mu$ of minimal degree,  so that $\mu$ is an inner node of $\ttt$ with $\deg(\mu)=\deg(\phi)$.


\subsection{Maclane--Vaqui\'e chains}\label{subsecMLV}

Consider a countable chain of valuations in $\ttt$: 
\begin{equation}\label{depthMLV}
v\ \stackrel{\phi_0,\ga_0}\lra\  \mu_0\ \stackrel{\phi_1,\ga_1}\lra\  \mu_1\ \stackrel{\phi_2,\ga_2}\lra\ \cdots
\ \stackrel{\phi_{n},\ga_{n}}\lra\ \mu_{n}\ \lra\ \cdots
\end{equation}
in which $\phi_0\in\kx$ is a monic polynomial of degree one, $\mu_0=[v;\,\phi_0,\ga_0]$ is a depth-zero valuation, and each other node is an augmentation  of the previous node, of one of the two types:\e

\emph{Ordinary augmentation}: \ $\mu_{n+1}=[\mu_n;\, \phi_{n+1},\ga_{n+1}]$, for some $\phi_{n+1}\in\kp(\mu_n)$.\e

\emph{Limit augmentation}:  \ $\mu_{n+1}=[\cc_n;\, \phi_{n+1},\ga_{n+1}]$,  for some $\phi_{n+1}\in\kpi(\cc_n)$, where $\cc_n$ is an essential continuous family  whose first valuation is $\mu_n$.\e

Therefore, $\phi_n$ is a key polynomial for $\mu_n$ of minimal degree and  $\deg(\mu_n)=\deg(\phi_n)$, for all $n\ge0$.\e

\defn
A chain of mixed augmentations as in (\ref{depthMLV}) is said to be  a \emph{\mlv (MLV) chain}  if every augmentation step satisfies:
\begin{itemize}
	\item If $\,\mu_n\to\mu_{n+1}\,$ is ordinary, then $\ \deg(\mu_n)<\deg(\mu_{n+1})$.
	\item If $\,\mu_n\to\mu_{n+1}\,$ is limit, then $\ \deg(\mu_n)=\deg(\cc_n)$ and $\ \phi_n\not \in\ty(\mu_n,\mu_{n+1})$. 
\end{itemize}\e

In this case, we have $\phi_n\not\in\ty(\mu_n,\mu_{n+1})$ for all $n$. As shown in \cite[Sec. 4.1]{MLV}, this implies that $\mu(\phi_n)=\ga_n$  and
$\g_{\mu_{n}}=\gen{\g_{\mu_{n-1}},\ga_{n}}$ for all $n$.

The following theorem is due to Maclane, for the discrete rank-one case  \cite{mcla}, and Vaqui\'e for the general case  \cite{Vaq}. Another proof may be found in  \cite[Thm. 4.3]{MLV}. 

\begin{theorem}\label{main}
Every node $\mu\in\ttt$ falls in one, and only one, of the following cases.  \e

(a) \ It is the last valuation of a finite MLV chain.
$$ \mu_0\ \stackrel{\phi_1,\ga_1}\lra\  \mu_1\ \stackrel{\phi_2,\ga_2}\lra\ \cdots\ \lra\ \mu_{r-1}\ \stackrel{\phi_{r},\ga_{r}}\lra\ \mu_{r}=\mu.$$

(b) \ It is the stable limit of an essential continuous family,  $\cc=\cfa$, whose first valuation $\mu_r$ falls in case (a):
$$ \mu_0\ \stackrel{\phi_1,\ga_1}\lra\ \mu_1\ \stackrel{\phi_2,\ga_2}\lra\   \cdots\ \lra\ \mu_{r-1}\ \stackrel{\phi_{r},\ga_{r}}\lra\ \mu_{r}\  \stackrel{\cfa}\lra\  \rho_{\cc}=\mu.
$$
Moreover, we may assume that $\deg(\mu_r)=\deg(\cc)$ and  $\phi_r\not\in\ty(\mu_r,\mu)$.\e

(c) \ It is the stable limit, $\mu=\lim_{n\in\N}\mu_n$, of an  infinite MLV chain.
\end{theorem}

The main advantage of MLV chains is that their nodes are essentially unique, so that we may read in them several data intrinsically associated to the valuation $\mu$.

For instance, the sequence $\left(\deg(\mu_n)\right)_{n\ge0}$ and the character ``ordinary" or ``limit" of each augmentation step $\mu_n\to\mu_{n+1}$, are intrinsic features of $\mu$ \cite[Sec. 4.3]{MLV}.  

Thus, we may define order preserving functions
$$
\dep,\ldp\colon\ttt\lra\N\infty,
$$
where $\dep(\mu)$ is the length of the MLV chain  underlying $\mu$, and
$\ldp(\mu)$ counts the number of limit augmentations in this MLV chain. 

It is easy to construct examples of valuations on $\kx$ of infinite depth. 
In the next section, we show the existence of valuations with infinite limit-depth too. Their construction is much more involved.

\section{A valuation with an infinite limit-depth}\label{secLimDepthInfinity}

In this section, we exhibit an example of a valuation with an infinite limit-depth, based on explicit constructions by Kuhlmann, of infinite  towers of field extensions with defect \cite{KuhlDefect}.

For a prime number $p$, let $\F$ be an algebraic closure of the prime field $\F_p$. For an indeterminate $t$, consider the fields of Laurent series, Newton-Puiseux series and Hahn series in $t$, respectively:
$$
\F((t))\subset K=\bigcup_{N\in\N} \F((t^{1/N}))\subset H=\F((t^\Q)).
$$
For a generalized power series $s=\sum_{q\in\Q}a_qt^q$, its support is a subset of $\Q$: 
$$\supp(s)=\{q\in\Q\mid a_q\ne0\}.$$

The Hahn field $H$ consists of all generalized power series with well-ordered support. The Newton-Puiseux field $K$ contains all series whose support is included in $\frac1N\Z_{\ge m}$ for some $N\in\N$, $m\in\Z$.

From now on, we denote by $\irr_K(b)$ the minimal polynomial over $K$ of any $b\in\kb$.

On these three fields we may consider the valuation $v$ defined as $$v(s)=\min(\supp(s)),$$
which clearly satisfies, $$v\left(\F((t))^*\right)=\Z, \qquad v(K^*)=v(H^*)=\Q.$$ The valued field $\left(\F((t)),v\right)$ is henselian, because it is the completion of the  discrete rank-one  valued field $\left(\F(t),v\right)$.
Since the extension $\F((t))\subset K$ is algebraic, the valued field $(K,v)$ is henselian too.

The Hahn field $H$ is algebraically closed. Thus, it contains an algebraic closure $\kb$ of $K$. The algebraic  generalized power series have been described by Kedlaya \cite{Ked0, Ked1}. Let us recall  \cite[Lem. 3]{Ked0}, which is essential for our purposes.

\begin{lemma}\label{kedlaya}
If $s\in H$ is algebraic over $K$, then it is contained in a tower of Artin-Schreier extensions of $K$. In particular, $s$ is separable over $K$ and $\deg_K s$ is a power of $p$.
\end{lemma}

Any $s\in H$ determines a valuation on $H[x]$ extending $v$:
$$
v_s\colon H[x]\lra\Q\infty,\qquad g\longmapsto v_s(g)=v(g(s)).
$$

We are interested in the valuation on $\kx$ obtained by restriction of $v_s$, which we still denote by the same symbol $v_s$. 
If $s$ is algebraic over $K$ and $f=\irr_K(s)\in\kx$, we have $v_s(f)=\infty$. 
Hence, $v_s$ cannot be extended to a valuation on $K(x)$.


On the other hand, suppose that $s=\sum_{q\in\Q}a_qt^q\in H$ is transcendental over $K$ and all its truncations 
$$
s_r=\sum_{q<r}a_qt^q, \qquad r\in\R,
$$
are algebraic over $K$ and have a bounded degree over $K$. Then, it is an easy exercise to check that $v_s$ falls in case (b) of Theorem \ref{main}.

Therefore, our example of a valuation with infinite limit-depth must be given by a transcendental $s\in H$, all whose truncations are algebraic over $K$ and have unbounded degree over $K$. In this case, $v_s$ will necessarily  fall in case (c) of Theorem \ref{main}. We want to find an example such that, moreover, all steps in the MLV chain of $v_s$ are limit augmentations. 

By Lemma \ref{kedlaya}, the truncations of $s$ must belong to some tower of Artin-Schreier extensions of $K$. Let us use a concrete tower constructed by Kuhlmann \cite[Ex. 3.14]{KuhlDefect}.

\subsection{A tower of Artin-Schreier extensions of $K$}
Let $\as(g)=g^p-g$ be the Artin-Screier operator on $\kx$. It is $\F_p$-linear and has kernel $\F_p$.

Let us start with the classical Abhyankar's example
$$
s_0=\sum_{i\ge1}t^{-1/p^i}\in H,
$$
which is a root of the polynomial $\varphi_0=\as(x)-t^{-1}\in \kx$. Since the denominators of the support of $s_0$ are unbounded, we have $s_0\not\in K$. Since the roots of $\varphi_0$ are $s_0+\ell$, for $\ell$ running on $\F_p$, the polynomial $\varphi_0$ has no roots in $K$. Hence, $\varphi_0$ is irreducible in $\kx$, because all irreducible polynomials in $\kx$ have degree a power of $p$.

Now, we iterate this construction to obtain a tower of Artin-Schreier extensions
$$
K\subsetneq K(s_0)\subsetneq K(s_1) \subsetneq\cdots \subsetneq K(s_n)\subsetneq \cdots
$$
where $s_n\in H$ is taken to be a root of $\varphi_n=\as(x)-s_{n-1}$. The above argument shows that $\varphi_n$ is irreducible in $K(s_{n-1})$ as long as $s_n\not \in K(s_{n-1})$, which is easy to check.   

From the algebraic relationship $\as(s_n)=s_{n-1}$ we may deduce  a concrete choice for all $s_n$:
$$
s_n=\sum_{j\ge n} \comb{j}{n}t^{-1/p^{j+1}},\quad\mbox{ for all }n\ge0,
$$
which follows from the well-known identity
$$
\comb{j+1}{n+1}=\comb{j}{n+1}+\comb{j}{n},\quad\mbox{ for all }j\ge n.
$$
In particular,
$$
\deg_K s_n= p^{n+1},\qquad v(s_n)=-1/p^{n+1},\quad\mbox{ for all }n\ge0.
$$

For all $n\ge0$, we have $\irr_K(s_n)=\as^{n+1}(x)-t^{-1}$, and the set of roots of this polynomial is 
\begin{equation}\label{ZAS}
Z\left(\irr_K(s_n)\right)=s_n+\op{Ker}(\as^{n+1})\subset s_n+\F.
\end{equation}
In particular, the support of all these conjugates of $s_n$ is contained in $(-1,0]$, and Krasner's constant of $s_n$ is zero:
\begin{equation}\label{krasner}
 \d(s_n)=\max\left\{v(s_n-\sg(s_n))\mid \sg\in\gal(\kb/K),\ \sg(s_n)\ne s_n\right\}=0.
\end{equation}

We are ready to define our transcendental $s\in H$ as:
$$
s=\sum_{n\ge0}t^ns_n.
$$

Let us introduce some useful notation to deal with the support of $s$ and its truncations. Consider the well-ordered set
$$
S=\left\{(n,i)\in\Z^2_{\lx}\mid \ 0\le n\le i,\quad p\nmid \comb{i}{n}\right\}.
$$
The support of $s$ is the image of the following order-preserving embedding
$$
\dta\colon S\hra \Q,\qquad (n,i)\longmapsto \dta(n,i)=n-\dfrac1{p^{i+1}}.
$$

The limit elements in $S$ are $(n,n)$ for $n\ge0$. These elements have no immediate predecessor in $S$. On the other hand, all elements in $S$ have an immediate successor:
$$
(n,i)\ \rightsquigarrow\ (n,i+m),
$$
where $m$ is the least natural number such that $p\nmid\comb{i+m}{n}$.

For all $(n,i)\in S$ we consider the truncations of $s$ determined by the rational numbers $\dta(n,i)$:
$$
\sni:=s_{\dta(n,i)}=\sum_{m=0}^{n-1}t^ms_m+t^n\sum_{j=n}^{i-1}\comb{j}{n}t^{-1/p^{j+1}}.
$$
For  the limit indices $(n,n)\in S$ the truncations are:
$$
s_{n,n}=\sum_{m=0}^{n-1}t^ms_m.
$$
Since $(0,0)=\min(S)$, the truncation $s_{0,0}=0$ is defined by an empty sum.

All truncations of $s$ are algebraic over $K$. Their degree is
$$
\deg_K \sni= p^n,\quad\mbox{ for all }(n,i)\in S,
$$
because $s_{n-1}$ has degree $p^n$, and all other summands have strictly smaller degree. For instance, the ``tail" \ $t^n\sum_{j=n}^{i-1}\comb{j}{n}t^{-1/p^{j+1}}$ belongs to $K$.

The unboundedness of the degrees of the truncations of $s$ is not sufficient to guarantee that $s$ is transcendental over $K$. To this end, we must analyze some more properties of these truncations. 

For any pair $(a,\dta)\in \kb\times \Q$, consider the ultrametric ball
$$
B=B_\dta(a)=\{b\in\kb\mid v(b-a)\ge\dta\}.
$$
We define the \emph{degree} of such a ball over $K$ as 
$$
\deg_K B=\min\{\deg_K b\mid b\in B\}.
$$

\begin{lemma}\label{minpair2}
For all $n\ge 1$, we have $\deg_K B_{n-1}\left(s_{n,n}\right)=p^n$. 
\end{lemma}

\begin{proof}
Denote $B=B_{n-1}\left(s_{n,n}\right)$.
From the computation in (\ref{krasner}), we deduce that Krasner's constant of $s_{n,n}$ is $\d(s_{n,n})=n-1$. 
Any $u\in B$  may be written as  
$$
u=s_{n,n}+\ell\, t^{n-1}+b, \qquad \ell\in\F,\quad b\in \kb,\quad v(b)>n-1.
$$
Let $z=s_{n,n}+\ell\, t^{n-1}$. Since $\ell\, t^{n-1}$ belongs to $K$, we have 
$$\deg_K z=p^n,\qquad \d(z)=n-1.$$ Since $v(u-z)>\d(z)$, we have $K(z)\subset K(u)$ by Krasner's lemma. Hence, $\deg_K u\ge p^n$. Since $B$ contains elements of degree $p^n$, we conclude that $\deg_K B=p^n$.  
\end{proof}

\begin{corollary}\label{strans}
The element $s\in H$ is transcendental over $K$. 
\end{corollary}

\begin{proof}
If $s$ were algebraic over $K$, it would belong to $B_{n-1}\left(s_{n,n}\right)$ for all $n$. This is impossible, because $\deg_K s$ would be unbounded, by Lemma \ref{minpair2}. 
\end{proof}

\subsection{A MLV chain of $v_s$ as a valuation on $\kb[x]$}

For all $(s,i)\in S$, we have
$$
v_s(x-\sni)=v(s-\sni)=n-\dfrac1{p^{i+1}}=\dta(n,i).
$$

Let $\vni=\om_{\sni,\dta(n,i)}$ be the depth zero valuation on $\kb[x]$ associated to the pair $(\sni,\dta(n,i))\in\kb\times\Q$; that is,
$$
\vni\left(\sum_{0\le\ell}a_\ell\left(x-\sni\right)^\ell\right)=\min_{0\le \ell}\left\{v_s\left(a_\ell\left(x-\sni\right)^\ell\right)\right\}=\min_{0\le \ell}\left\{v\left(a_\ell\right)+\ell\,\dta(n,i)\right\}.
$$

\begin{lemma}\label{vivs}
For all $(n,i),(m,j)\in S$ we have
$\vni(x-s_{m,j})=\min\{\dta(n,i),\dta(m,j)\}$. In particular,
$\vni<v_s$ for all $(n,i)\in S$.
\end{lemma}

\begin{proof}
 The computation of $\vni(x-s_{m,j})$ follows immediately form the definition of $\vni$. The inequality $\vni\le v_s$ follows from the comparison of the action of both valuation on  $(x-\sni)$-expansions. Finally, if we take $\dta(n,i)<\dta(m,j)$, we get
 $$
 \vni(x-s_{m,j})=\dta(n,i)<\dta(m,j)=v_s(x-s_{m,j}). $$
This shows that $\vni<v_s$. 
\end{proof}

\begin{lemma}\label{vsb}
The family $\cc=\left(\vni\right)_{(n,i)\in S}$ is a totally ordered family of valuations on $\kb[x]$ of stable degree one, admitting $v_s$ as its stable limit. 
\end{lemma}

\begin{proof}
Let us see that $\cc$ is a totally ordered family of valuations. More precisely, 
$$
(n,i)<(m,j)\ \imp\ \dta(n,i)<\dta(m,j)\ \imp\ \vni<v_{m,j}<v_s. 
$$
 Indeed, this follows from (\ref{balls}) because $v(\sni-s_{m,j})=v(\sni-s)=\dta(n,i)$. 

Clearly, $\cc$ contains no maximal element, and all valuations in $\cc$ have degree one.  Let us show that all polynomials $x-a\in\kb[x]$ are $\cc$-stable, and the stable value coincides with $v_s(x-a)=v(s-a)$.

Since $s$ is transcendental over $K$, we have $s\ne a$ and $q=v(s-a)$ belongs to $\Q$. For all $(n,i)\in S$ such that $\dta(n,i)>q$ we have
$$
\vni(x-a)=\min\{v(a-\sni),\dta(n,i)\}=\min\{q,\dta(n,i)\}=q=v_s(x-a).
$$
This ends the proof of the lemma.
\end{proof}\e

Therefore, $v_s$ falls in case (b) of Theorem \ref{main}, as a valuation on $\kb[x]$. 
A MLV chain of $v_s$ is, for instance,
$$
v_{0,0}\stackrel{\cc}\lra v_s=\lim(\cc).
$$

In order to obtain a MLV chain of $v_s$ as a valuation of $\kx$, we need to ``descend" this result to $\kx$. 
In this regard, we borrow some ideas of \cite{Vaq3}.

\subsection{A MLV chain of $v_s$ as a valuation on $\kx$}
We say that  $(a,\dta)\in \kb\times \Q$ is a  \emph{minimal pair} if $\deg_K B_\dta(a)=\deg_K a$. This concept was introduced in \cite{APZ}. By equation (\ref{balls}), for all $b\in\kb$ we have$$\om_{a,\dta}=\om_{b,\dta} \ \sii\ b\in B_\dta(a).$$However, only the minimal pairs $(a,\dta)$ of this ball contain all essential information about the valuation on $\kx$ that we obtain by restriction of $\om_{a,\dta}$. 

\begin{lemma}\cite[Prop. 3.3]{Vaq3}\label{minpair}
For $(a,\dta)\in\kb\times\Q$, let $\mu$ be the valuation on $\kx$ obtained by restriction of the valuation $\om=\om_{a,\dta}$ on $\kb[x]$. 
Then, for all $g\in\kx$, $\inn_\mu g$ is a unit in $\ggm$ if and only if   $\inn_{\om} g$ is a unit in $\gg_{\om}$.
\end{lemma}

The following result was originally proved in \cite{PP}; another proof can be found in \cite[Thm. 1.1]{N2019}.  

\begin{lemma}\label{minp}
For a minimal pair $(a,\dta)\in\kb\times\Q$, let $\mu$ be the valuation on $\kx$ obtained by restriction of the valuation $\om_{a,\dta}$ on $\kb[x]$. 
Then, $\irr_K(a)$ is a key polynomial for $\mu$, of minimal degree. 
\end{lemma}

We need a last auxiliary result.

\begin{lemma}\label{minpair3}
For all $(n,i)\in S$ the pair $\left(\sni,\dta(n,i)\right)$ is minimal. 
\end{lemma}

\begin{proof}
All  $\left(s_{0,i},\dta(0,i)\right)$ are minimal pairs, because $\deg_K s_{0,i}=1$. For $n>0$, denote $B_{n,i}=B_{\dta(n,i)}\left(\sni\right)$. Since $B_{n,i}\subset B_{n-1}\left(s_{n,n}\right)$, Lemma \ref{minpair2} shows that 
$$
\deg_K B_{n,i}\,\ge\,\deg_K B_{n-1}\left(s_{n,n}\right)=p^n.
$$
Since the center $s_{n,i}$ of the ball $B_{n,i}$ has $\deg_K \sni=p^n$, we deduce $\deg_K B_{n,i}=p^n$. Thus, $\left(\sni,\dta(n,i)\right)$ is a minimal pair.
\end{proof}\e

\noindent{\bf Notation. }Let us denote the restriction of $\vni$ to $\kx$ by
$$
\rni=\left(\vni\right)_{\mid \kx}.
$$
Moreover, for the limit indices $(n,n)$, $n\ge0$, we denote:
$$
\mu_n=\rho_{n,n},\qquad \phi_n=\irr_K\left(s_{n,n}\right),\qquad \ga_n=v_s\left(\phi_n\right). 
$$
By Lemmas \ref{vivs} and \ref{vsb}, the set of all valuations $\left(\rni\right)_{(n,i)\in S}$ is totally ordered, and $\rni<v_s$ for all $(n,i)$. 

\begin{proposition}\label{mainvs}
For all $n\ge0$, the set $\cc_n=\left(\rni\right)_{(n,i)\in S}$ is an essential continuous family of stable degree $p^n$.
Moreover, the polynomial $\phi_{n+1}$ belongs to $\kpi\left(\cc_n\right)$ and $\mu_{n+1}=[\cc_n;\,\phi_{n+1},\ga_{n+1}]$. 
\end{proposition}

\begin{proof}
Let us fix some $n\ge0$. By Lemmas \ref{minp} and \ref{minpair3}, all valuations in $\cc_n$ have degree $p^n$. Hence, $\cc_n$ is a totally ordered family of stable degree $p^n$. 

Let us show that all monic $g\in \kx$ with $\deg(g)\le p^n$ are $\cc_n$-stable. 
Let $u\in\kb$ be a root of $g$. By Lemma \ref{minpair2}, $u\not\in B_n\left(s_{n+1,n+1}\right)$, so that $v\left(s_{n+1,n+1}-u\right)<n$. Since $v\left(s-s_{n+1,n+1}\right)=\dta(n+1,n+1)>n$, we deduce that $v(s-u)<n$.

Therefore, we may find $j\ge n$ such that 
$$
v(s-u)<n-\dfrac1{p^{j+1}}
$$
for all roots $u$ of $g$. As we showed along the proof of Lemma \ref{vsb}, this implies $$\vni(x-u)=v_s(x-u)\quad\mbox{  for all }(n,i)\ge (n,j)$$ simultaneously for all roots $u$ of $g$. Therefore, $\rni(g)=v_s(g)$ for all $(n,i)\ge (n,j)$ and $g$ is $\cc_n$-stable.

Now, let us show that $\phi_{n+1}$ is $\cc_n$-unstable. For all $i\ge n$, we have 
$$v\left(s_{n+1,n+1}-\sni\right)=\dta(n,i)=\vni(x-\sni).$$ 

By \cite[Prop. 6.3]{KP}, $x-s_{n+1,n+1}$ is a key polynomial for $\vni$; thus, $\inn_{\vni} \left(x-s_{n+1,n+1}\right)$ is not a unit in the graded algebra $\gg_{\vni}$. Hence, $\inn_{\vni}\phi_{n+1}$ is not a unit in $\gg_{\vni}$ and  Lemma \ref{minpair} shows that 
$\inn_{\rni}\phi_{n+1}$ is not a unit in $\gg_{\rni}$. Since this holds for all $i$, Lemma \ref{stable=unit} shows that $\phi_{n+1}$ is $\cc_n$-unstable.

Since the irreducible polynomials in $\kx$ have degree a power of $p$ (Lemma \ref{kedlaya}), $\phi_{n+1}$ is an $\cc_n$-unstable polynomial of minimal degree. Therefore, $\cc_n$ is an essential continuous family and $\phi_{n+1}\in\kpi(\cc_n)$.

Since $\phi_{n+1}$ is $\cc_n$-unstable,  $\rni\left(\phi_{n+1}\right)<v_s\left(\phi_{n+1}\right)=\ga_{n+1}$ for all $i$. Thus, it makes sense to consider the limit augmentation $\mu=[\cc_n;\,\phi_{n+1},\ga_{n+1}]$. Let us show that $\mu=\mu_{n+1}$ by comparing their action on $\phi_{n+1}$-expansions. For all $g=\sum_{0\le \ell}a_\ell\phi_{n+1}^\ell$,
\begin{equation}\label{bothmu}
\mu_{n+1}(g)=\min_{0\le\ell}\left\{\mu_{n+1}\left(a_\ell\phi_{n+1}^\ell\right)\right\},\qquad\mu(g)=\min_{0\le\ell}\left\{\mu\left(a_\ell\phi_{n+1}^\ell\right)\right\}. 
\end{equation}

Since $\deg(a_\ell)<p^{n+1}=\deg\left(\phi_{n+1}\right)$, all these coefficients $a_\ell$ are  $\cc_n$-stable. Hence, $\rni\left(a_\ell\right)=v_s\left(a_\ell\right)$ for all $(n,i)$ sufficiently large.  Since $\rni<\mu_{n+1}<v_s$, we deduce  
 $$\mu\left(a_\ell\right)=\rho_{\cc_n}\left(a_\ell\right)=\rni\left(a_\ell\right)=\mu_{n+1}\left(a_\ell\right)=v_s\left(a_\ell\right).$$ 

Finally, for all $i\ge n+1$, we have $v\left(s_{n+1,i}-s_{n+1,n+1}\right)=\dta(n+1,n+1)$, so that
$$
v_{n+1,n+1}\left(x-s_{n+1,n+1}\right)=\dta(n+1,n+1)=v_{n+1,i}\left(x-s_{n+1,n+1}\right).
$$ 
By (\ref{ZAS}), for all the other roots $u$ of $\phi_{n+1}$, the support of $u$ is contained in $(-1,n]$. Thus, for all $i\ge n+1$ we get
$$
v_{n+1,n+1}\left(x-u\right)=v\left(s_{n+1,n+1}-u\right)=v\left(s_{n+1,i}-u\right)=v_{n+1,i}\left(x-u\right).
$$ 
Since $\mu_{n+1}=\rho_{n+1,n+1}<\rho_{n+1,i}<v_s$,  \cite[Cor. 2.5]{MLV} implies
$$
\mu_{n+1}\left(\phi_{n+1}\right)=\rho_{n+1,i}\left(\phi_{n+1}\right)=v_s\left(\phi_{n+1}\right)=\ga_{n+1}=\mu\left(\phi_{n+1}\right).
$$

By (\ref{bothmu}), we deduce that $\mu=\mu_{n+1}$.
\end{proof}\e

Therefore, we get a countable chain of limit augmentations
$$
\mu_0\ \stackrel{\phi_1,\ga_1}\lra\  \mu_1\ \stackrel{\phi_2,\ga_2}\lra\ \cdots
\ \lra\ \mu_{n-1} 
\ \stackrel{\phi_{n},\ga_{n}}\lra\ \mu_{n} \ \lra\ \cdots
$$
which is an MLV chain. Indeed, the MLV condition amounts to
$$
\phi_n\not\in \ty(\mu_n,\mu_{n+1})\quad \mbox{ for all }n\ge0.
$$
This means $\mu_n(\phi_n)=\mu_{n+1}(\phi_n)$ for all $n$. Since $\mu_n<\mu_{n+1}<v_s$, the desired equality  follows from $\mu_n(\phi_n)=\ga_n=v_s(\phi_n)$. 

Finally, the family $\left(\mu_n\right)_{n\in\N}$ has stable limit $v_s$. Indeed, for all nonzero $f\in\kx$, there exists $n\in\N$ such that $\deg(f)<p^n=\deg(\mu_n)$. Let $\ty(\mu_n,v_s)=[\varphi]_{\mu_n}$. Since $\deg(f)<\deg(\varphi)$, we have $\varphi\nmid_{\mu_n}f$ and this implies $\mu_n(f)=v_s(f)$ by Lemma \ref{propertiesTMN}.

As a consequence, $v_s$ has infinite limit-depth.



\end{document}